%&amstex
\input amsppt.sty
\magnification=\magstep1
\hsize=30truecc
\vsize=20truecm
\baselineskip=24truept
\NoBlackBoxes
\pageno=1
\TagsOnRight
\nologo
\def\Z{\Bbb Z}

\def\N{\Bbb N}

\def\al{\alpha}
\def\l{\left}
\def\r{\right}
\def\bg{\bigg}
\def\({\bg(}
\def\[{\bg[}
\def\){\bg)}
\def\]{\bg]}
\def\t{\text}
\def\f{\frac}

\def\mo{\roman{mod}}

\def\em{\emptyset}
\def\se {\subseteq}

\def\sm{\setminus}

\def\bi{\binom}
\def\eq{\equiv}
\def\cs{\cdots}
\def\ls{\leqslant}
\def\gs{\geqslant}
\def\mo{\roman{mod}}
\def\ord{\roman{ord}}

\def\Proof{\noindent{\it Proof}}
\def\Def{\medskip\noindent{\it Definition}\ }

\topmatter \hbox{A talk given at the Institute of Math. Science,
Nanjing Univ. on Oct. 8, 2004.}
\medskip
\title {Groups and Combinatorial Number Theory}\endtitle
\author Zhi-Wei Sun\endauthor
\affil Department of Mathematics\\Nanjing University\\Nanjing 210093, P. R. China
\\{\it E-mail:} {\tt zwsun\@nju.edu.cn}
\\Homepage: {\tt http://pweb.nju.edu.cn/zwsun}
\endaffil
\abstract
    In this talk we introduce several topics in combinatorial number theory
which are related to groups; the topics include
combinatorial aspects of covers of groups by cosets,
and also restricted sumsets and zero-sum problems on abelian groups.
A survey of known results and open
problems on the topics is given in a popular way.
\endabstract
\thanks 2000 {\it Mathematics Subject Classifications}:\,Primary 20D60;
Secondary 05A05, 11B25, 11B75, 11P70.
\newline\indent
All papers of the author mentioned in this survey are available from his homepage
{\tt http://pweb.nju.edu.cn/zwsun}.\endthanks
\endtopmatter
\document
\hsize=30truecc
\baselineskip=24truept

\heading{1. Nontrivial problems and results on cyclic groups}\endheading

Any infinite cyclic group is isomorphic to the additive group $\Z$ of all integers.
Subgroups of $\Z$ different from $\{0\}$ are those $n\Z=\{nq:\,q\in\Z\}$ with $n\in\Z^+=\{1,2,3,\ldots\}$.
Any cyclic group of order $n$ is isomorphic to the additive group $\Z/n\Z$ of residue classes modulo $n$.
A coset of the subgroup $n\Z$ of $\Z$ has the form
$$a+n\Z=\{a+nq:\, q\in\Z\}=\{x\in\Z:\ x\eq a\ (\mo\ n)\}$$
which is called a {\it residue class} with {\it modulus} $n$ or an arithmetic sequence with common difference $n$.
For convenience we also write $a(n)$ or $a(\mo\ n)$ for $a+n\Z$, thus $0(1)=\Z$ and $1(2)$ is the set of odd integers.

We can decompose the group $\Z$ into $n$ cosets of $n\Z$, namely
$$\{r(n)\}_{r=0}^{n-1}=\{0(n),\ 1(n),\ \ldots,\ n-1(n)\}$$
is a partition of $\Z$ (i.e., a disjoint cover of $\Z$). For the index of the subgroup $n\Z$ of $\Z$,
we clearly have $[\Z:n\Z]=|\Z/n\Z|=n$.

Since $0(2^n)$ is a disjoint union of the residue classes $2^n(2^{n+1})$ and $0(2^{n+1})$, the systems
$$\align&A_1=\{1(2),0(2)\},\ A_2=\{1(2),2(4),0(4)\},\ A_3=\{1(2),2(4),4(8),0(8)\},
\\&\qquad \cdots\cdots,\  A_k=\{1(2),2(2^2),\ldots,2^{k-1}(2^{k}), 0(2^{k})\},\ \cdots\cdots
\endalign$$
are disjoint covers of $\Z$.

The concept of cover of $\Z$ was first introduced by P. Erd\H os in the early 1930s. He noted that
$\{0(2),\ 0(3),\ 1(4),\ 5(6),\ 7(12)\}$
is a cover of $\Z$ with the moduli $2,3,4,6,12$ distinct.

Soon after his invention of the concept of cover of $\Z$, Erd\H os
made the following conjecture: {\it If $A=\{a_s(n_s)\}_{s=1}^k \ (k>1)$
is a system of residue classes with the moduli
$n_1,\ldots,n_k$ distinct, then it cannot be a disjoint cover of $\Z$}.

\proclaim{Theorem 1.1} Let $A=\{a_s(n_s)\}_{s=1}^k$.

{\rm (i) (H. Davenport, L. Mirsky, D. Newman and R. Rad\'o)} If $A$
is a disjoint cover of $\Z$ with $1<n_1\ls n_2\ls\cdots\ls n_{k-1}\ls n_k$,
then we must have $n_{k-1}=n_k$.

{\rm (ii) [Z. W. Sun, Chinese Quart. J. Math. 1991]} Let $n_0$ be a positive period of the function
$w_A(x)=|\{1\ls s\ls k:\, x\in a_s(n_s)\}|$. For any positive integer $d$
with $d\nmid n_0$ and $I(d)=\{1\ls s\ls k:\,d\mid n_s\}\not=\em$, we have
$$|I(d)|\gs|\{a_s\ \mo\ d:\, s\in I(d)\}|
\gs\min\Sb 0\ls s\ls k\\d\nmid n_s\endSb\f d{\gcd(d,n_s)}\gs p(d),$$
where $p(d)$ is the least prime divisor of $d$.
\endproclaim

\noindent{\it Proof of part $(i)$}. Without loss of generality we let $0\ls a_s<n_s$ ($1\ls s\ls k$).
For $|z|<1$ we have
$$\sum^k_{s=1}\f{z^{a_s}}{1-z^{n_s}}=\sum^k_{s=1}\sum^{\infty}_{q=0}z^{a_s+qn_s}=\sum^{\infty}_{n=0}
z^n=\f1{1-z}.$$
If $n_{k-1}<n_k$ then
$$\infty=\lim\Sb z\to e^{2\pi i/n_k}\\ |z|<1\endSb\f {z^{a_k}}{1-z^{n_k}}=\lim\Sb z\to e^{2\pi i/n_k}\\ |z|<1
\endSb\l(\f1{1-z}-\sum^{k-1}_{s=1}\f{z^{a_s}}{1-z^{n_s}}\r)<\infty,$$
a contradiction! \qed

Part (ii) in the case $n_0=1$ and $d=n_k$ yields the Davenport-Mirsky-Newman-Rad\'o result,
a further extension of part (ii)
was given by Z. W. Sun [Math. Res. Lett. 11(2004)] and [J. Number Theory, to appear].

Recall that
$$A_k=\{1(2),2(2^2),\ldots,2^{k-1}(2^{k}), 0(2^{k})\}$$
is a disjoint cover of $\Z$. Thus the system
$\{1(2),2(2^2),\ldots,2^{k-1}(2^{k})\}$
covers $1,\ldots,2^k-1$ but does not cover any multiple of $2^k$.
In 1965 P. Erd\H os made the following conjecture.

\proclaim{Erd\H os' Conjecture} $A=\{a_s(n_s)\}_{s=1}^k$ forms a cover of $\Z$
 if it covers those integers from $1$ to $2^k$.
\endproclaim

In 1969--1970 R. B. Crittenden and
C. L. Vanden Eynden [Bull. Amer. Math. Soc. 1969; Proc. Amer. Math. Soc. 1970]
supplied a long and awkward proof
of the Erd\H os conjecture for $k\gs20$,
which involves some deep results concerning the distribution of primes.

The following result is stronger than Erd\H os' conjecture.
\proclaim{Theorem 1.2 {\rm [Z. W. Sun, Acta Arith. 72(1995),
Trans. Amer. Math. Soc. 348(1996)]}}
Let $A=\{a_s(n_s)\}_{s=1}^k$ be a finite system of residue classes,
and let $m_1,\ldots,m_k$ be integers relatively prime to
$n_1,\ldots,n_k$ respectively. Then system $A$
forms an $m$-cover of $\Z$ (i.e., $A$ covers every integer at least $m$ times)
if it covers $|S|$ consecutive
integers at least $m$ times, where
$$S=\bg\{\bg\{\sum_{s\in I}\f{m_s}{n_s}\bg\}:\, I\se\{1,\ldots,k\}\bg\}.$$
$($As usual the fractional part of a real number $x$
is denoted by $\{x\}.)$
\endproclaim
\noindent{\tt Proof of Theorem 1.2 in the case $m=1$}.  For any integer $x$, clearly
$$\align&x\ \t{is covered by}\ A
\\\iff& e^{2\pi i(a_s-x)m_s/n_s}=1\ \t{for some}\ s=1,\ldots,k
\\\iff&\prod_{s=1}^k\l(1-e^{2\pi i(a_s-x)m_s/n_s}\r)=0
\\\iff&\sum_{I\se\{1,\ldots,k\}}(-1)^{|I|}
e^{2\pi i\sum_{s\in I}a_sm_s/n_s}\cdot e^{-2\pi ix\sum_{s\in I}m_s/n_s}=0
\\\iff&\sum_{\theta\in S}e^{-2\pi ix\theta}z_{\theta}=0,
\endalign$$
where
$$z_{\theta}=\sum\Sb I\se\{1,\ldots,k\}\\\{\sum_{s\in I}m_s/n_s\}=\theta\endSb
(-1)^{|I|}e^{2\pi i\sum_{s\in I}a_sm_s/n_s}.$$

Suppose that $A$ covers $|S|$ consecutive integers
$a,a+1,\ldots,a+|S|-1$ where $a\in\Z$. By the above,
$$\sum_{\theta\in S}(e^{-2\pi i\theta})^r(e^{-2\pi ia\theta}z_{\theta})=0$$
for $r=0,1,\ldots,|S|-1$. As the determinant
$\|(e^{-2\pi i\theta})^r\|_{0\ls r<|S|,\,\theta\in S}$
is of Vandermonde's type and hence nonzero, by Cramer's rule
we have $z_{\theta}=0$ for all $\theta\in S$.
Therefore $\sum_{\theta\in S}e^{-2\pi ix\theta}z_{\theta}=0$ for all $x\in\Z$,
i.e., any $x\in\Z$ is covered by $A$. This proves the theorem in the case $m=1$.
 \qed

The following theorem shows that disjoint covers of $\Z$ are related to unit fractions,
actually further results were obtained by Z. W. Sun.

\proclaim{Theorem 1.3 {\rm (Z. W. Sun [Acta Arith. 1995; Trans. Amer. Math. Soc. 1996])}} Let
$A=\{a_s(n_s)\}_{s=1}^k$ be a disjoint cover of $\Z$.

{\rm (i)} If $\em \not=J\subset\{1,\ldots,k\}$, then there exists an $I\se\{1,\ldots,k\}$ with $I\not=J$
such that $\sum_{s\in I}1/n_s=\sum_{s\in J}1/n_s$.

{\rm (ii)} For any $1\ls t\ls k$ and $r\in\{0,1,\ldots,n_t-1\}$, there
is an $I\se\{1,\ldots,k\}\sm\{t\}$ such that $\sum_{s\in I}1/n_s=r/n_t$.
\endproclaim
\Proof. Let $N=[n_1,\ldots,n_k]$ be the least common multiple of $n_1,\ldots,n_k$. Then
$$\prod_{s=1}^k\l(1-z^{N/n_s}e^{2\pi ia_s/n_s}\r)=1-z^N$$
because each $N$th root of unity is a single zero of the left hand side.
Thus
$$\sum_{I\se\{1,\ldots,k\}}(-1)^{|I|}z^{\sum_{s\in I}N/n_s}e^{2\pi i\sum_{s\in I}a_s/n_s}=1-z^N.$$
Comparing the degrees of both sides we obtain the well-known equality $\sum_{s=1}^k1/n_s=1$.
As $\em \not=J\subset\{1,\ldots,k\}$, $0<\sum_{s\in J}N/n_s<N$ and hence
$$\sum\Sb I\se\{1,\ldots,k\}\\\sum_{s\in I}1/n_s=\sum_{s\in J}1/n_s\endSb
(-1)^{|I|}e^{2\pi i\sum_{s\in I}a_s/n_s}=0$$
which implies that $\sum_{s\in I}1/n_s=\sum_{s\in J}1/n_s$ for some $I\se\{1,\ldots,k\}$ with $I\not=J$.

Now fix $1\ls t\ls k$. Observe that
$$\prod^k\Sb s=1\\s\not=t\endSb\l(1-z^{N/n_s}e^{2\pi ia_s/n_s}\r)=\f{1-z^N}{1-z^{N/n_t}e^{2\pi ia_t/n_t}}
=\sum_{r=0}^{n_t-1}z^{Nr/n_t}e^{2\pi ia_tr/n_t}.$$
Thus, for any $r=0,1,\ldots,n_t-1$ we have
$$\sum\Sb I\se\{1,\ldots,k\}\sm\{t\}\\\sum_{s\in I}1/n_s=r/n_t\endSb(-1)^{|I|}e^{2\pi i\sum_{s\in I}a_s/n_s}
=e^{2\pi ia_tr/n_t}\not=0$$
and hence $\sum_{s\in I}1/n_s=r/n_t$ for some $I\se\{1,\ldots,k\}\sm\{t\}$. \qed

We mention that covers of $\Z$ by residue classes have many surprising applications.
For example, on the basis of Cohen and Selfridge's work,
Z. W. Sun [Proc. Amer. Math. Soc. 2000] showed that {\it  if
$$x\eq 47867742232066880047611079\ (\mo\ M)$$
then $x$ is not of the form $\pm p^a\pm q^b$, where $p,q$ are
primes and $a,b$ are nonnegative integers, and $M$ is a 29-digit number given by
$$\align&\prod_{p\ls 19}p\times31\times37\times41\times61\times73\times97
\times109\times151\times241\times257\times331
\\&\qquad\ \ =66483084961588510124010691590.
\endalign$$}

If $A=\{a_1<\cdots<a_k\}$ and $B=\{b_1<\cdots<b_l\}$ are finite subsets of $\Z$, then
clearly the sumset $A+B=\{a+b:\,a\in A\ \&\ b\in B\}$ contains at least the following $k+l-1$ elements:
$$a_1+b_1<a_2+b_1<\cdots<a_k+b_1<a_k+b_2<\cdots<a_k+b_l.$$
However, the following result for
cyclic groups of prime orders is nontrivial and very useful.
\proclaim{Theorem 1.4 {\rm (Cauchy-Davenport Theorem)}} Let $A$ and $B$ be nonempty
subsets of $\Z_p=\Z/p\Z$ where $p$ is a prime. Then we have
$$|A+B|\gs\min\{p,|A|+|B|-1\}.$$
\endproclaim

 In 1964 P. Erd\H os and Heilbronn posed the following conjecture for cyclic groups
 of prime orders.

\proclaim{Erd\H os-Heilbronn Conjecture} Let $A$ be a nonempty
subset of $\Z_p=\Z/p\Z$ where $p$ is a prime. Then we have
$$|2^{\wedge}A|\gs\min\{p,2|A|-3\}$$
where $2^{\wedge}A=\{a+b:\,a,b\in A\ \&\ a\not=b\}$.
\endproclaim

This conjecture remained open until it was confirmed by Dias da Silva and
Y. Hamidoune [Bull. London. Math. Soc. 1994]
thirty years later, with the help of the representation theory of groups.

\proclaim{Theorem 1.5 {\rm (The da Silva--Hamidoune Theorem)}}
  Let $p$ be a prime and $\em\not=A\se\Z_p$. Then we have
$$|n^{\wedge}A| \ge \min\{p,\ n|A|-n^2+1\},$$
where $n^{\wedge}A$ denotes the set of all sums of $n$ distinct elements
of $A$.
\endproclaim

If $p$ is a prime, $A\se\Z_p$ and $|A|>\sqrt{4p-7}$, then
by the da Silva--Hamidoune theorem, any element of $\Z_p$
can be written as a sum of $\lfloor|A|/2\rfloor$ distinct elements of $A$.

In 1995--1996  Alon, Nathanson and Ruzsa [Amer. Math. Monthly 1995, J. Number Theory 1996]
developed a polynomial method rooted in [Alon and Tarsi, Combinatorica 1989]
to prove the Erd\H os-Heilbronn conjecture and some similar results. The method
turns out to be very powerful and has many applications
in number theory and combinatorics.

An extension of Theorem 1.5 appeared in
Q. H. Hou and Z. W. Sun [Acta Arith. 2002].
H. Pan and Z. W. Sun [J. Combin. Theory Ser. A 2002] obtained a general result
on sumsets with polynomial restrictions which includes the Cauchy-Davenport theorem as
a special case.

 Suppose that
 $$\{a_1,\cs,a_n\},\ \{b_1,\cs,b_n\}\ \t{and}\ \{a_1+b_1,\cs,a_n+b_n\}$$
are complete systems of residues modulo $n$.
Let
$$\sigma=0+1+\cs+(n-1)=\f{n(n-1)}2.$$
As
$$\sum_{i=1}^n(a_i+b_i)=\sum_{i=1}^na_i+\sum_{i=1}^nb_i,$$
we have
$\sigma\eq\sigma+\sigma\ (\mo\ n)$ and hence $2\nmid n$.

 In 2001, Dasgupta, K\'arolyi, Serra and
 Szegedy [Israel J. Math. 2001] confirmed a conjecture of H. S. Snevily
 for cyclic groups.

 \proclaim{Theorem 1.6 {\rm (Dasgupta-K\'arolyi-Serra-Szegedy Theorem)}}
\ Let $G$ be an additive cyclic
group with $|G|$ odd. Let $A$ and $B$ be subsets of $G$ with
cardinality $n>0$. Then there is a numbering $\{a_i\}_{i=1}^n$ of
the elements of $A$ and a numbering $\{b_i\}_{i=1}^n$ of the
elements of $B$ such that $a_1+b_1,\cs,a_n+b_n$ are pairwise
distinct.
 \endproclaim
 \Proof. As $2^{\varphi(|G|)}\eq1\ (\mo\ |G|)$ (where $\varphi$ is Euler's totient function),
the {\bf multiplicative} group of the finite field $F$ with order $2^{\varphi(|G|)}$
has a cyclic subgroup isomorphic to $G$.
Thus we can view $G$ as a subgroup of the multiplicative group $F^*=F\sm\{0\}$.

 Write $A=\{a_1,\ldots,a_n\}$ and $B=\{b_1,\ldots,b_n\}$. We want to show that
 there is a $\sigma\in S_n$ such that
 $a_{\sigma(i)}b_i\not=a_{\sigma(j)}b_j$ whenever $1\ls i<j\ls n$.
 In other words,
 $$c=\sum_{\sigma\in S_n}\prod_{1\ls i<j\ls n}\l(a_{\sigma(j)}b_j-a_{\sigma(i)}b_i\r)\not=0.$$
 In fact,
 $$\align c=&\sum_{\sigma\in S_n}\|a_{\sigma(j)}^{i-1}b_j^{i-1}\|_{1\ls i,j\ls n}\ \ (\t{Vandermonde})
 \\=&\sum_{\sigma\in S_n}\sum_{\tau\in S_n}\t{sign}(\tau)\prod_{j=1}^na_{\sigma(j)}^{\tau(j)-1}b_j^{\tau(j)-1}
 \\=&\sum_{\tau\in S_n}\t{sign}(\tau)
 \prod_{j=1}^nb_j^{\tau(j)-1}\sum_{\sigma\in S_n}\prod_{j=1}^na_{\sigma(j)}^{\tau(j)-1}
\\=&\sum_{\tau\in S_n}\t{sign}(\tau)\prod_{j=1}^nb_j^{\tau(j)-1}\sum_{\sigma\in S_n}\t{sign}(\sigma\tau^{-1})
\prod_{i=1}^na_{\sigma\tau^{-1}(i)}^{i-1}\ (\t{as}\ -1=1\ \t{in}\ F)
\\=&\|b_j^{i-1}\|_{1\ls i,j\ls n}\times\|a_j^{i-1}\|_{1\ls i,j\ls n}=\prod_{1\ls i<j\ls n}(a_j-a_i)(b_j-b_i)\not=0.
\endalign$$
 This concludes the proof. \qed

 The following conjecture remains unsolved.

\proclaim{Snevily's Conjecture} Let $a_1,\ldots,a_k\in\Z$, and let $n$ be a positive integer greater than $k$.
Then there is a permutation $\sigma\in S_k$ such that
all the $i+a_{\sigma(i)}\ (i=1,\ldots,k)$ modulo $n$ are distinct.
\endproclaim

\heading{2. Nontrivial Problems and Results on Abelian Groups}\endheading

Let $G$ be an additive abelian group of order $n$,
and let $b_1,\cs,b_n\in G$. If both
$\{a_i\}_{i=1}^n$ and $\{a_i+b_i\}_{i=1}^n$
are numberings of the elements of $G$, then
$\sum_{i=1}^n(a_i+b_i)=\sum_{i=1}^na_i$ and hence $b_1+\cs+b_n=0$.
In 1952 M. Hall [Proc. Amer. Math. Soc.] obtained the converse.

\proclaim{Theorem 2.1 {\rm (Hall's Theorem)}} Let
$G=\{a_1,\cs,a_n\}$ be an additive abelian group, and let
$b_1,\cs,b_n$ be any elements of $G$ with $b_1+\cs+b_n=0$.
Then there exists a permutation $\sigma\in S_n$
such that $a_{\sigma(1)}+b_1,\cs,a_{\sigma(n)}+b_n$
are distinct.
\endproclaim

Hall's proof is highly technical.

In 1999 H. S. Snevily made the following general conjecture:

{\it  Let $G$ be any additive abelian group with $|G|$ odd. Let $A$ and $B$ be subsets of $G$ with
cardinality $n>0$. Then there is a numbering $\{a_i\}_{i=1}^n$ of
the elements of $A$ and a numbering $\{b_i\}_{i=1}^n$ of the
elements of $B$ such that $a_1+b_1,\cs,a_n+b_n$ are pairwise
distinct.}

 The proof of the following result in this direction involves linear algebra, field theory and
 Dirichlet's unit theorem in algebraic number theory

\proclaim{Theorem 2.2 {\rm (Z. W. Sun [J. Combin. Theory Ser. A, 2003])}}
 Let $G$ be an additive abelian group whose finite
 subgroups are all cyclic. Let $A_1,\cs,A_n$ be finite
 subsets of $G$ with cardinality $k>m(n-1)$
 (where $m$ is a positive integer), and let $b_1,\cs,b_n$
 be elements of $G$.

 {\rm (i)} If $b_1,\cs,b_n$ are distinct,
 then there are at least $(k-1)n-m\bi n2+1$
 multi-sets $\{a_1,\cs,a_n\}$ such that $a_i\in A_i$ for $i=1,\cs,n$
 and all the $ma_i+b_i$ are distinct.

 {\rm (ii)} The sets
  $$\{\{a_1,\cs,a_n\}\colon a_i\in A_i,\ a_i\not=a_j\
 \t{and}\ ma_i+b_i\not= ma_j+b_j\ \t{if}\ i\not=j\}$$
 and $$\{\{a_1,\cs,a_n\}\colon a_i\in A_i,\ ma_i\not=ma_j\
 \t{and}\ a_i+b_i\not= a_j+b_j\ \t{if}\ i\not=j\}$$
 have more than $(k-1)n-(m+1)\bi n2\gs(m-1)\bi n2$ elements,
 provided that $b_1,\cs,b_n$ are distinct and of odd order,
 or they have finite order and
 $n!$ cannot be written in the form $\sum_{p\in P}px_p$ where all the $x_p$
 are nonnegative integers and $P$ is the set of primes dividing
 one of the orders of $b_1,\cs,b_n$.
 \endproclaim

 In the 1960's M. Knerser obtained the following remarkable theorem on abelian groups.

\proclaim{Theorem 2.3 {\rm (Kneser's Theorem)}} Let $G$ be an additive abelian group.
Let $A$ and $B$ be finite nonempty subsets of $G$, and let
$H=H(A+B)$ be the stablizer $\{g\in G:\, g+A+B=A+B\}$. If $|A+B|\ls|A|+|B|-1$, then
$$|A+B|=|A+H|+|B+H|-|H|.$$
\endproclaim

The following consequence is an extension of the Cauchy-Davenport theorem.

\proclaim{Corollary 2.1} Let $G$ be an additive abelian group. Let $p(G)$
be the least order of a nonzero element of $G$, or $p(G)=+\infty$ if $G$ is torsion-free.
Then, for any finite nonempty subsets $A$ and $B$ of $G$, we have
$$|A+B|\gs\min\{p(G),|A|+|B|-1\}.$$
\endproclaim
\Proof. Suppose that $|A+B|<|A|+|B|-1$. Then $H=H(A+B)\not=\{0\}$ by Kneser's theorem.
Therefore $|H|\gs p(G)$ and hence
$$|A+B|=|A+H|+|B+H|-|H|\gs|A+H|\gs|H|\gs p(G).$$
We are done. \qed

Quite recently G. K\'arolyi was able to extend
the Erd\H os-Heilbronn conjecture to any abelian groups.

\proclaim{Theorem 2.4 {\rm (G. K\'arolyi [Israel J. Math. 2004])}} Let $G$ be an additive
abelian group. Then, for any finite nonempty subset $A$ of $G$, we have
$$|2^{\wedge}A|\gs\min\{p(G),2|A|-3\}.$$
\endproclaim

The characteristic function of a residue class is a periodic arithmatical map.
Dirichlet characters are also periodic functions. If an element $a$
in an additive abelian group $G$ has order $n$, then the map $\psi:\Z\to G$
given by $\psi(x)=xa$
is periodic mod $n$.

\proclaim{Theorem 2.5 {\rm (Z. W. Sun, 2004)}}
Let $G$ be any additive abelian group, and let $\psi_1,\ldots,\psi_k$
be maps from $\Z$ to $G$ with periods $n_1,\ldots,n_k\in\Z^+$ respectively. Then
the function $\psi=\psi_1+\cdots+\psi_k$ is constant if
 $\psi(x)$ equals a constant for $|T|\ls n_1+\cdots+n_k-k+1$ consecutive integers $x$, where
 $$T=\bigcup_{s=1}^k\l\{\f r{n_s}:\, r=0,1,\ldots,n_s-1\r\}.$$
 \endproclaim

 The proof of Theorem 2.5 involves linear recurrences and algebraic integers.
 \proclaim{Corollary 2.2 {\rm (Z. W. Sun [Math. Res. Lett. 11(2004)])}} The system
 $A=\{a_s(\mo\ n_s)\}_{s=1}^k$ covers every integer exactly $m$ times
 if it covers $|T|$ consecutive integers exactly $m$ times, where $T$ is as in
 Theorem 2.5.
 \endproclaim

 In 1966 J. Mycielski [Fund. Math.] posed an interesting conjecture
on disjoint covers (i.e. partitions) of abelian groups.
Before stating the conjecture we give a definition first.

\Def\ 2.1. The Mycielski function $f:\Z^+=\{1,2,\ldots\}\to \N=\{0,1,2,\ldots\}$
is given by
$$f(n)=\sum_{p\in P(n)}\ord_p(n)(p-1),$$
where $P(n)$ denotes the set of prime divisors of $n$ and $\ord_p(n)$
represents largest integer $\al$ such that $p^{\al}\mid n$.
In other words, $f(\prod_{t=1}^rp_t^{\al_t})=\sum_{t=1}^r\al_t(p_t-1)$ where $p_1,\ldots,p_r$
are distinct primes.
\medskip

\proclaim{Mycielski's Conjecture} Let $G$ be an abelian group, and $\{a_sG_s\}_{s=1}^k$
be a disjoint cover of $G$ by left cosets of subgroups.
Then $k\gs 1+f([G:G_t])$ for each $t=1,\ldots,k$.
(It is known that $[G:G_t]<\infty$ for all $t=1,\ldots,k$.)
\endproclaim

Mycielski's conjecture was first confirmed by \v S. Zn\'am [Colloq. Math., 1966] in the case $G=\Z$.

\proclaim{Theorem 2.6 {\rm (G. Lettl and Z. W. Sun, 2004)}}
 Let $\Cal A=\{a_sG_s\}_{s=1}^k$
be a cover of an abelian group $G$ by left cosets of subgroups. Suppose that
$\Cal A$ covers all the elements of $G$ at least $m$ times with the coset $a_tG_t$
irredundant. Then $[G:G_t]\ls 2^{k-m}$ and furthermore
$k\gs m+f([G:G_t])$.
\endproclaim

In the case $m=1$ and $G_t=\{e\}$, this confirms a conjecture of W. D. Gao
and A. Geroldinger [European J. Combin. 2003].
The proof of Theorem 2.6 involves algebraic number theory and characters of abelian groups.

\proclaim{Conjecture {\rm (Z. W. Sun, 2004)}} Let $\Cal A=\{a_sG_s\}_{s=1}^k$
be a finite system of left cosets of subgroups of an abelian group $G$. Suppose that
$\Cal A$ covers all the elements of $G$ at least $m$ times but none of its proper subsystems does.
Then we have $k\gs m+f(N)$ where $N$ is the least common multiple of the indices $[G:G_1],\ldots,[G:G_k]$.
\endproclaim

In 1961 P. Erd\H os, A. Ginzburg and A. Ziv [Bull. Research Council. Israel] established the following
celebrated theorem which initiated the study of zero-sums.

\proclaim{Theorem 2.7 {\rm (The EGZ Theorem)}} Let $G$ be any additive abelian group of
order $n$. For any given $c_1,\cs,c_{2n-1}\in G$, there is an $I\se\{1,\ldots,2n-1\}$
with $|I|=n$ such that $\sum_{s\in I}c_s=0$.
\endproclaim

In 2003 Z. W. Sun connected the EGZ theorem with covers of $\Z$.

\proclaim{Theorem 2.8 {\rm (Z. W. Sun [Electron. Res. Announc. AMS, 2003])}}
Let $A=\{a_s(\mo\ n_s)\}_{s=1}^k$
and suppose that $|\{1\ls s\ls k:\,x\eq a_s\ (\mo\ n_s)\}|\in\{2q-1,2q\}$
for all $x\in\Z$, where $q$ is a prime power. Let $G$ be an additive abelian group
of order $q$. Then, for any $c_1,\ldots,c_k\in G$, there exists an $I\se\{1,\ldots,k\}$
such that $\sum_{s\in I}1/n_s=q$ and $\sum_{s\in I}c_s=0$.
\endproclaim

\Def\ 2.2. The Davenport constant $D(G)$ of a finite abelian group $G$
(written additively) is defined as the smallest positive
integer $k$ such that any sequence $\{c_s\}_{s=1}^k$ (repetition
allowed) of elements of $G$ has a nonempty subsequence
$c_{i_1},\cs,c_{i_l}\ (i_1<\cs<i_l)$ with zero-sum (i.e.
$c_{i_1}+\cs+c_{i_l}=0$).
\medskip

For any abelian group $G$ of order $n$ we clearly have $D(G)\ls n$.
In fact, if $c_1,\ldots,c_n\in G$, then the partial sums
$$s_0=0,\ s_1=a_1,\ s_2=a_1+a_2,\ \ldots,\ s_n=a_1+\cdots+a_n$$
cannot be distinct since $n+1>|G|$,
so there are $0\ls i<j\ls n$ such that $s_i=s_j$, i.e. $a_{i+1}+\cdots+a_j=0$.

In 1966 Davenport showed that if $K$ is
an algebraic number field with ideal class group $G$, then $D(G)$
is the maximal number of prime ideals (counting multiplicity) in
the decomposition of an irreducible integer in $K$.

In 1969 J. Olson
[J. Number Theory] used the knowledge of group rings to show that
the Davenport constant of an abelian $p$-group $G\cong
\Z_{p^{h_1}}\oplus\cs\oplus\Z_{p^{h_l}}$ is
$1+\sum_{t=1}^l(p^{h_t}-1).$

 In 1994 W. R. Alford, A. Granville and C. Pomerance [Ann. Math.]
 employed an upper bound for the Davenport constant of
 the unit group of the ring $\Z_n$ to prove that
 there are infinitely many Carmichael numbers which are those composites
 $m$ such that $a^{m-1}\eq1\ (\mo\ m)$ for any $a\in\Z$ with $(a,m)=1$.

 The following well-known conjecture is still open, it is known to be true for $k=1,2$.

 \proclaim{Olson's Conjecture} Let $k$ and $n$ be positive integers. Then
 $D(\Z_n^k)=1+k(n-1)$
 where $\Z_n^k$ is the direct sum of $k$ copies of $\Z_n$.
 \endproclaim

 \heading {3. Nontrivial Problems and Results on General Groups}\endheading

 Let $G$ be a group and $G_1,\cs,G_k$ be subgroups of $G$.
Let $a_1,\cs,a_k\in G$. If the system
$\Cal A=\{a_iG_i\}_{i=1}^k$ of left cosets
covers all the elements of $G$ at least $m$ times
but none of its proper subsystems does, then
all the indices $[G:G_i]$ are known to be finite.

\proclaim{Theorem 3.1} Let $\Cal A=\{a_iG_i\}_{i=1}^k$
be a finite system of left cosets in a group $G$ where $G_1,\ldots,G_k$
are subgroups of $G$. Suppose that $\Cal A$ forms a minimal cover $G$
(i.e. $\Cal A$ covers all the elements of $G$
but none of its proper systems does).

{\rm (i) (B. H. Neumann [Publ. Math. Debrecen, 1954])} There is a constant $c_k$ depending only on $k$
such that $[G:G_i]\ls c_k$ for all $i=1,\ldots,k$.

{\rm (ii) (M. J. Tomkinson [Comm. Algebra, 1987]) We have
$[G:\bigcap_{i=1}^kG_i]\ls k!$ where the upper bound $k!$ is best possible.
\endproclaim
\Proof. We prove (ii) by induction. (Part (ii) is stronger than part (i).)

 We want to show that
 $$\[\bigcap_{i\in I}G_i:\bigcap_{i=1}^kG_i\]\ls (k-|I|)!\tag $*_I$ $$
 for all $I\se\{1,\ldots,k\}$,
 where $\bigcap_{i\in\em}G_i$ is regarded as $G$.

 Clearly $(*_I)$ holds for $I=\{1,\ldots,k\}$.

  Now let $I\subset\{1,\ldots,k\}$ and assume $(*_J)$ for all $J\se\{1,\ldots,k\}$ with $|J|>|I|$.
 Since $\{a_iG_i\}_{i\in I}$ is not a cover of $G$, there is an $a\in G$ not covered by $\{a_iG_i\}_{i\in I}$.
 Clearly $a(\bigcap_{i\in I} G_i)$ is disjoint from the union $\bigcup_{i\in I}a_iG_i$
 and hence contained in $\bigcup_{j\not\in I}a_jG_j$. Thus
 $$a\(\bigcap_{i\in I}G_i\)=\bigcup\Sb j\not\in I\\a_jG_j\cap a(\bigcap_{i\in I}G_i)\not=\em\endSb
 \(a_jG_j\cap a\(\bigcap_{i\in I}G_i\)\)$$
 and hence
 $$\[\bigcap_{i\in I}G_i:H\]\ls\sum_{j\not\in I}\[G_j\cap\bigcap_{i\in I}G_i:H\]
 \ls\sum_{j\not\in I}(k-(|I|+1))!=(k-|I|)!$$
 where $H=\bigcap_{i=1}^kG_i$.
 This concludes the induction proof. \qed

 \medskip
 \Def\ 3.1. Let $H$ be a subnormal subgroup of a group $G$ with finite index,
and
$$H_0=H\subset H_1\subset\cs\subset H_n=G$$
be a composition series from $H$ to $G$
(i.e. $H_i$ is maximal normal in $H_{i+1}$ for each $0\ls i<n$).
If the length $n$ is zero (i.e. $H=G$), then
we set $d(G,H)=0$, otherwise we put
$$d(G,H)=\sum_{i=0}^{n-1}([H_{i+1}:H_i]-1).$$
\medskip

Let $H$ be a subnormal subgroup of a group $G$ with $[G:H]<\infty$.
By the Jordan--H\"older theorem,
$d(G,H)$ does not depend on the choice of the
composition series from $H$ to $G$.
Clearly $d(G,H)=0$ if and only if $H=G$. If $K$ is a subnormal subgroup
of $H$ with $[H:K]<\infty$, then
$$d(G,H)+d(H,K)=d(G,K).$$
When $H$ is normal in $G$,
the `distance' $d(G,H)$ was first introduced by I. Korec [Fund. Math. 1974].
The current general notion is due to Z. W. Sun [Fund. Math. 1990].
Z. W. Sun [Fund. Math. 1990]
showed that
$$[G:H]-1\gs d(G,H)\gs f([G:H])\gs \log_2[G:H]$$
where $f$ is the Mycielski function.
Moreover, Sun [European J. Combin. 2001] noted that
$d(G,H)=f([G:H])$ if and only if $G/H_G$ is solvable
where $H_G=\bigcap_{g\in G}gHg^{-1}$ is the largest normal
subgroup of $G$ contained in $H$.

In 1968 \v S. Znam [Coll. Math. Soc. J\'anos Bolyai] made the following further conjecture:
{\it If $A=\{a_s(n_s)\}_{s=1}^k$ is a disjoint cover of $\Z$ then
$$k\gs1+f(N_A)\ \ \t{and hence}\ N_A\ls2^{k-1},$$
where $N_A=[n_1,\ldots,n_k]=[\Z:\bigcap_{s=1}^kn_s\Z]$.}

In 1974 I. Korec [Fund. Math.] confirmed Zna\'am's conjecture and Mycielski's conjecture
by proving the following deep result: {\it Let $\{a_iG_i\}^k_{i=1}$ be
a partition of a group into left cosets of
normal subgroups. Then
$k\gs1+f([G:\bigcap^k_{i=1}G_i]).$}

Here is a further extension of Korec's result.

\proclaim{Theorem 3.2 {\rm (Z. W. Sun [European J. Combin. 22(2001)])}}
Let $G$ be a group and $\{a_iG_i\}_{i=1}^k$
cover each elements of $G$ exactly $m$ times, where $G_1,\ldots,G_k$ are
subnormal subgroups of $G$. Then
$$k\gs m+d\bg(G,\bigcap_{i=1}^kG_i\bg),$$
where the lower bound can be attained.
Moreover, for any subgroup $K$ of $G$ not contained in all the $G_i$ we have
$$|\{1\ls i\ls k:K\not\se G_i\}|\gs 1+d\(K,K\cap\bigcap_{i=1}^kG_i\).$$
\endproclaim

\proclaim{Corollary 3.1 {\rm (Z. W. Sun [Fund. Math. 1990])}}
Let $H$ be a subnormal subgroup of a group $G$ with $[G:H]<\infty$.
Then
$$[G:H]\gs1+d(G,H_G)\gs 1+f([G:H_G])\ \ \t{and hence}\ [G:H_G]\ls 2^{[G:H]-1}.$$
\endproclaim
\Proof. Let $\{Ha_i\}_{i=1}^k$ be a right coset decomposition of $G$ where $k=[G:H]$.
Then $\{a_iG_i\}_{i=1}^k$ is a disjoint cover of $G$ where all
the $G_i=a_i^{-1}Ha_i$ are subnormal in $G$. Observe that
$$\bigcap_{i=1}^kG_i=\bigcap_{i=1}^k\bigcap_{h\in H}a_i^{-1}h^{-1}Hha_i
=\bigcap_{g\in G}g^{-1}Hg=H_G.$$
So the desired result follows from Theorem 3.2. \qed

\proclaim{Theorem 3.3} {\rm (i) (Berger-Felzenbaum-Fraenkel, 1988, Coll. Math.)}
If  $\{a_iG_i\}^k_{i=1}$ is a  disjoint cover of a finite solvable group
$G$, then $k\gs 1+f([G:G_i])$ for $i=1,\cs,k$.

{\rm (ii) [Z. W. Sun, European J. Combin. 2001]} Let $G$ be a group and
$\{a_iG_i\}_{i=1}^k$ be a finite system of left cosets
which covers each elements of $G$ exactly $m$ times.
For any $i=1,\ls,k$, whenever $G/(G_i)_G$ is solvable
we have  $k\gs m+f([G:G_i])$ and hence $[G:G_i]\ls2^{k-m}$.
\endproclaim

Z. W. Sun [European J. Combin. 2001] suggested the following further conjecture.

\proclaim{Conjecture 3.1 {\rm (Z. W. Sun, 2001)}}
Let $a_1G_1,\ldots,a_kG_k$ be left cosets of a group $G$
such that $\{a_iG_i\}_{i=1}^k$ covers each elements of $G$ exactly $m$ times and that
all the $G/(G_i)_G$ are solvable. Then $k\gs m+f(N)$ where $N$ is the least common multiple
of the indices $[G:G_1],\cs,[G:G_k]$.
\endproclaim

If $\{x_1,\ldots,x_k\}$ is a maximal subset of a group $G$ with $x_ix_j\not=x_jx_i$ for all $1\ls i<j\ls k$,
then $\{C_G(x_i)\}_{i=1}^k$ is a minimal cover of $G$ with $\bigcap_{i=1}^kC_G(x_i)=Z(G)$
(Tomkinson, Comm. Algebra, 1987)
and $|G/Z(G)|\ls c^k$ for some absolute constant (L. Pyber, J. London Math. Soc., 1987).

\proclaim{Conjecture 3.2 {\rm (Z. W. Sun, 1996)}}
Let $\{G_i\}_{i=1}^k$ be a minimal cover of a group $G$ by subnormal subgroups.
Write $[G:\bigcap_{i=1}^kG_i]=\prod_{t=1}^rp_t^{\al_t}$, where $p_1,\ldots,p_r$ are distinct primes
and $\al_1,\ldots,\al_r$ are positive integers.
Then we have
$$k\gs1+\sum_{t=1}^r(\al_t-1)p_t.$$
\endproclaim

Up to now, no counterexample to this conjecture has been found.

The following conjecture extends a conjecture of P. Erd\H os.

\proclaim{The Herzog-Sch\"onheim Conjecture {\rm ([Canad. Math. Bull. 1974])}}
\ Let $\Cal A=\{a_iG_i\}_{i=1}^k\ (k>1)$
be a partition (i.e. disjoint cover) of a group $G$
into left cosets of subgroups $G_1,\cs,G_k$. Then
the indices $n_1=[G:G_1],\cs,n_k=[G:G_k]$
 cannot be pairwise distinct.
\endproclaim

M. A. Berger, A. Felzenbaum and A. S. Fraenkel [1986, Canad. Math. Bull.; 1987, Fund. Math.]
showed the conjecture for finite nilpotent groups and supersolvable
groups. A quite recent progress was made by the speaker.

\proclaim{Theorem 3.4 {\rm (Z. W. Sun [J. Algebra, 2004])}}
 Let $G$ be a group, and
 $\Cal A=\{a_iG_i\}_{i=1}^k$ $(k>1)$
 be a system of left cosets of subnormal subgroups.
 Suppose that $\Cal A$ covers each $x\in G$ the same number of times, and
$$n_1=[G:G_1]\ls\cs\ls n_k=[G:G_k].$$
Then the indices $n_1,\cs,n_k$ cannot be
distinct. Moreover, if each index occurs in $n_1,\cs,n_k$
at most $M$ times, then
$$\log n_1\ls\f{e^{\gamma}}{\log 2}M\log^2M+O(M\log M\log\log M)$$
where $\gamma=0.577\cs $ is the Euler constant and the $O$-constant is absolute.
\endproclaim

The above theorem also
answers a question analogous to a famous problem of Erd\H os negatively.
Theorem 3.4 was established by a combined use of tools from
 group theory and number theory.

One of the key lemmas is the following one which is the main reason
why covers involving subnormal subgroups are better behaved than general covers.

 \proclaim{Lemma 3.1 {\rm (Z. W. Sun [European J. Combin. 2001])}} Let $G$ be a group,
 and let $P(n)$ denote the set of prime divisors of a positive integer $n$.

 {\rm (i)} If $G_1,\ldots,G_k$ are subnormal subgroups of $G$ with finite index, then
 $$\[G:\bigcap_{i=1}^kG_i\]\ \big|\ \prod_{i=1}^k[G:G_i]
 \ \t{and hence}\ P\(\[G:\bigcap_{i=1}^kG_i\]\)=\bigcup_{i=1}^kP([G:G_i]).$$

 {\rm (ii)} Let $H$ be a subnormal subgroup of $G$ with finite index. Then
 $$P(|G/H_G|)=P([G:H]).$$
\endproclaim

We mention that part (ii) is a consequence of the first part, and the word ``{\it subnormal}"
cannot be removed from part (i).

Here is another useful lemma.

\proclaim{Lemma 3.2 {\rm (Z. W. Sun [J. Algebra, 2004])}} Let $G$ be a group
and $H$ its subgroup with finite index $N$. Let $a_1,\ldots,a_k\in G$,
and let $G_1,\ldots,G_k$
be subnormal subgroups of $G$ containing $H$. Then $\bigcup_{i=1}^ka_iG_i$
contains at least $|\bigcup_{i=1}0(n_i)\cap\{0,1,\ldots,N-1\}|$ left cosets of $H$,
where $n_i=[G:G_i]$.
\endproclaim

Finally we pose an interesting unsolved conjecture.

\proclaim{Conjecture 3.3 {\rm (Z. W. Sun)}} Let $G$ be a group,
and $a_1G_1,\ldots,a_kG_k$ be pairwise disjoint left cosets of $G$
with all the indices $[G:G_i]$ finite. Then, for some $1\ls i<j\ls k$ we have
$\gcd([G:G_i],[G:G_j])\gs k$.
\endproclaim

This conjecture is open even in the special case $G=\Z$.

\enddocument